\documentclass[reqno]{amsart} %

\usepackage[marginratio=1:1,height=660pt,width=460pt,tmargin=60pt]{geometry}  

\usepackage{amssymb}
\usepackage{hyperref}
\usepackage[utf8]{inputenc} 
\usepackage{amssymb,amsmath,amsthm}
\usepackage{cite}
\usepackage{hyperref} 
\usepackage{verbatim} 
\usepackage{graphicx,subfig,epstopdf}  

\newcommand{\ipa}{\ensuremath{[Re(\alpha)]+1}}

\AtBeginDocument{{\noindent\small
\sc{Bulletin of Mathematical Analysis and Applications}\newline
ISSN: 1821-1291, URL: http://www.bmathaa.org\newline
Volume 6 Issue 4 (2014), Pages 1-15}

\thanks{\copyright 2014 Universiteti i Prishtin\"es, Prishtin\"e, Kosov\"e.}
\vspace{9mm}}

\begin{document}

\begin{center} \large{\textbf{A NEW APPROACH TO GENERALIZED FRACTIONAL DERIVATIVES}}

\vspace{.8cm}

Udita N. Katugampola

\vspace{.3cm}

\small{Department of Mathematics, University of Delaware, Newark, DE 19716, USA} \vspace{-.3cm}

\end{center}

\title[ Generalized Fractional Derivatives]{}


\author[U.N. Katugampola]{}

\address{Udita N. Katugampola \newline
 Department of Mathematics,
 University of Delaware,
 Newark, DE 19716, USA}
 \email{uditanalin@yahoo.com}


\thanks{Submitted December 2, 2013. Published October 15, 2014 \hfill Email:uditanalin@yahoo.com}



\begin{abstract}
The author \mbox{(Appl. Math. Comput. 218(3):860-865, 2011)} introduced a new fractional integral operator given by, 
\[
\big({}^\rho \mathcal{I}^\alpha_{a+}f\big)(x) = \frac{\rho^{1- \alpha }}{\Gamma({\alpha})} \int^x_a \frac{\tau^{\rho-1} f(\tau) }{(x^\rho - \tau^\rho)^{1-\alpha}}\, d\tau, 
\] 
which generalizes the well-known Riemann-Liouville and the Hadamard fractional integrals. In this paper we present a new fractional derivative which generalizes the familiar Riemann-Liouville and the Hadamard fractional derivatives to a single form. We also obtain two representations of the generalized derivative in question. An example is given to illustrate the results.
\end{abstract}

\maketitle

\numberwithin{equation}{section}
\newtheorem{remark}{Remark}

\newtheorem{theorem}{Theorem}[section]
\newtheorem{corollary}[theorem]{Corollary}
\newtheorem{proposition}[theorem]{Proposition}
\newtheorem{lemma}[theorem]{Lemma}
\newtheorem{conclusion}[theorem]{Conclusion} 
\newtheorem{question}[theorem]{Question}
\theoremstyle{definition}
\newtheorem{definition}[theorem]{Definition}
\newtheorem{example}[theorem]{Example}
\newtheorem{problem}[theorem]{Problem}


\section{Introduction}\label{sec:1}

In recent years, the Fractional Calculus (FC) draws increasing attention due to its applications in many fields. The history of the theory goes back to seventeenth century, when in 1695 the derivative of order $\alpha = \frac{1}{2}$ was described by Leibnitz in his letter to L'Hospital \cite{letter1,letter2,letter3}. Since then, the new theory turned out to be very attractive to mathematicians as well as physicists, biologists, engineers and economists. The first application of fractional calculus was due to Abel in his solution to the Tautocrone problem \cite{Abel}. It also has applications in biophysics, quantum mechanics, wave theory, polymers, continuum  mechanics, Lie theory, field theory, spectroscopy and in group theory, among other applications \cite{Herr,Hilf,hilferthreefold,what}. In \cite{key-9}, Samko et al. provide an encyclopedic treatment of the subject. Various type of fractional derivatives were studied: Riemann-Liouville, Caputo, Hadamard, Erd\'{e}lyi-Kober, Gr\"{u}nwald-Letnikov, Marchaud and Riesz are just a few to name \cite{key-9,key-8,key-2,key-11,key-12}.  

In fractional calculus, the fractional derivatives are defined via fractional integrals \cite{key-9,key-8}. According to the literature, the Riemann-Liouville fractional derivative (RLFD), hence the Riemann-Liouville fractional integral plays a major role in FC \cite{key-9}. The Caputo fractional derivative has also been defined via a modified Riemann-Liouville fractional integral \cite{key-8}. Butzer et al. investigate properties of the Hadamard fractional integral and the derivative in \cite{key-9,key-8,key-2,key-1,key-5,key-6,key-3,key-7}. In \cite{key-6}, they also obtained the Mellin transforms of the Hadamard fractional integral and diffferential operators and in \cite{Pooseh}, Pooseh et al. obtained expansion formulas of the Hadamard operators in terms of integer order derivatives. Many other interesting properties of those operators and others are summarized in \cite{key-9} and \cite{key-8} and the references therein. 

In \cite{key-0}, the author introduced a new fractional integral, which generalizes the Riemann-Liouville and the Hadamard integrals into a single form. For further properties such as expansion formulas, variational calculus applications, control theoretical applications, convexity and integral inequalities and Hermite-Hadamard type inequalities of this new operator and similar operators, for example, see \cite{Herr,Pooseh,u-1,u-2,u-3,u-4,u-5,u-6,u-7,u-8,u-9}. In the present work, we shall introduce a new fractional derivative, which generalizes the two derivatives in question. 

The paper is organized as follows. In the next section, we give definitions and some properties of the fractional integrals and fractional derivatives of various types. More detailed explanation can be found in the book by Samko et al. \cite{key-9} and the references therein.

\section{Definitions}

We shall start this section with some historical remarks and definitions to refresh our memories about some of the remarkable milestones in the theory of fractional calculus. As is well known nowadays, the first documented note about a fractional derivative was found in 1695 in the letters of Leibnitz to L'Hospital \cite{letter1,letter2,letter3,Ross1}. In 1819, Lacroix obtained the well-known $\frac{1}{2}$ derivative of $x$ \cite{lacroix}, using inductive arguments, to be $d^{\frac{1}{2}}/dx^{\frac{1}{2}} = 2\sqrt{x/\pi}$, long before the Riemann-Liouville fractional derivative surfaced into the realm of fractional calculus. The idea of a derivative that is not of an order of a positive integer was introduced by Liouville in 1832 \cite{liouville1,Lutzen,Osler1}, in a manner that would generalize the relation $D^ne^{\alpha x} = \alpha^n e^{\alpha x}$ to any complex number $n$. Liouville then used Fourier theory to extend his $\alpha$-derivative to any function $f(z)$ expanded in a Fourier series \cite{Osler1}. In 1888, Nekrassov \cite{Nekrassov,Osler1,Osler2}, generalizing the Cauchy's integral formula, 
\[
    \frac{d^n f(z)}{dz^n}=\frac{n!}{2\pi i}\oint_C \frac{f(\zeta)}{(\zeta -z)^{n+1}}\,d\zeta,
\]
where $C$ is a closed contour surrounding the point $z$ and enclosing a region of analyticity of $f$, came up with a fractional derivative which can be showed to be equal to Riemann-Liouville derivative under certain conditions \cite[p.54-55]{key-12}.

The \emph{Riemann-Liouville fractional integrals} $I^{\alpha}_{a+}f$ and $I^{\alpha}_{b-}f$ of order $\alpha \in \mathbb{C},$ $(Re(\alpha)>0)$ are defined by \cite{key-9,liouville1,key-8,Riemann1,riesz1},
\begin{equation}
(I^\alpha_{a+} f)(x) = \frac{1}{\Gamma(\alpha)}\int_a^x (x - \tau)^{\alpha -1} f(\tau) d\tau  \quad ; x > a,
\label{eq:lRLI}
\end{equation}
and 
\begin{equation}
(I^\alpha_{b-} f)(x) = \frac{1}{\Gamma(\alpha)}\int_x^b (\tau - x)^{\alpha -1} f(\tau) d\tau  \quad ; x < b,
\label{eq:rRLI}
\end{equation}
respectively. Here $\Gamma(\cdot)$ is the Gamma function. These integrals are called the \emph{left-sided} and \emph{right-sided}  fractional integrals, respectively. When $\alpha = n \in \mathbb{N},$ the integrals (\ref{eq:lRLI}) and (\ref{eq:rRLI}) coincide with the $n$-fold integrals \cite[chap.2]{key-8}. 
The corresponding \emph{Riemann-Liouville fractional derivatives} $D^{\alpha}_{a+}f$ and $D^{\alpha}_{b-}f$ of order $\alpha \in \mathbb{C},\, Re(\alpha)\geq 0$ are defined by \cite{key-9},
\begin{equation}
(D^\alpha_{a+} f)(x)=\bigg(\frac{d}{dx}\bigg)^n \, \Big(I^{n-\alpha}_{a+} f\Big)(x), \quad  x>a,
\label{eq:lRLD}
\end{equation}
and
\begin{equation}
(D^\alpha_{b-} f)(x)=\bigg(-\frac{d}{dx}\bigg)^n \, \Big(I^{n-\alpha}_{b-} f\Big)(x),  \quad  x<b,
\label{eq:rRLD}
\end{equation}
respectively, where $n = \ipa $. A word about notations is necessary here. We sometimes use the ceiling function, $\lceil \cdot \rceil$ to denote the same quantity $[\cdot]+1$, when there is no room for confusion. For simplicity, from this point onwards, we consider only the \emph{left-sided} integrals and derivatives, except in a few occasions. The interested reader may find more detailed information about \emph{right-sided} integrals and derivatives in the references, for example in \cite{key-9,key-8}.

One of the disadvantages of RLFD is that it is not consistant with the physical initial and boundary conditions when it comes to initial or boundary value problems. To overcome this difficulty, M. Caputo coined a variation of RLFD, now known in the litureture as \emph{Caputo} or \textit{Dzherbashyan-Caputo} fractional derivative given by \cite{caputo1,dzher1,key-8},
\begin{equation*}
({}^cD^\alpha_{a+} f)(x)=\frac{1}{\Gamma(n-\alpha)}\int_a^x (x - \tau)^{n-\alpha -1} f^{(n)}(\tau) d\tau, \quad  n-1<\alpha\leq n.
\label{eq:caputo}
\end{equation*}

Some historical notes about the derivative in question can be found in \cite[p. 18-21]{mainardi1}. The interested reader may also find an extended reference lists about Caputo derivative, for example, in the book by Mainardi \cite{mainardi1}. An application oriented treatment of the Caputo derivative is given in the book by Diethelm \cite{diethlem1}. As pointed out by Hilfer \cite{hilferthreefold}, the Caputo derivative was originally introduced by Liouville \cite[p.10]{Liouville2} though it did not take much attention until Caputo brought the idea back to life in his celebrated paper \emph{``Linear models of dissipation whose Q is almost frequency independent, Part II''}\cite{caputo}.

An interpolation between the two derivatives mentioned above are defined by the \emph{Hilfer fractional derivative} of order $\alpha$ and type $\beta$ given by \cite[p.113]{Hilf}\cite[p.11]{mainardi1},
\[
    {}_0D^{\alpha,\beta}_t := I^{\beta(1-\alpha)}_{0+} \circ D^1 \circ I^{(1-\beta)(1-\alpha)}_{0+}, \;\; 0 <\alpha, \beta \leq 1. 
\]
The Riemann-Liouville derivative of order $\alpha$ corresponds to the type $\beta=0$, while the Caputo derivative to the type $\beta=1$.

The \emph{Weyl-Riesz} fractional integration operator of a periodic function $f$ takes the form \cite{Momani1},
\[
   W^\alpha = \frac{1}{2\pi}\int_0^{2\pi}\Psi^\alpha(\tau)f(x-\tau)\,d\tau,
\]
where
\[
 \Psi^\alpha(\tau)=2\sum_{n=1}^{\infty}\frac{\cos n\tau}{n^\alpha}, \; \mbox{and} \; 0 \leq \tau \leq 2\pi.
\]
Further properties of this operator can be found in \cite{key-9,Samko2} and the references therein. Jumarie proposed a simple modification to the Riemann-Liouville derivative given by \cite{Jumarie1,Jumarie2,atan}, 
\[
   D^\alpha_x f(x)=\frac{1}{\Gamma(n-\alpha)} \frac{d^n}{dx^n}\int_0^x (x - \tau)^{n-\alpha -1} \big[f(\tau)-f(0)\big] d\tau, \quad  n-1<\alpha \leq n.
\]
with the property that the derivative of a constant function being equal to zero, a property important in Applied Mathematics and Engineering applications. 

Another fractional integral that is important in potential theory is the Riesz fractional integral given by \cite{Agrawal1},
\begin{equation}
   {}^RI_t^\alpha f(t) = \frac{1}{2\Gamma(\alpha)}\int_a^b |t - \tau|^{\alpha -1} f(\tau) \,d\tau, \; \; \; \alpha > 0.
	\label{RI}
\end{equation}
Two other variations of this integral, one with a factor of $1/\cos(\pi\alpha/2)$ and another with a factor of $2$ can be found in \cite{key-9} and \cite{key-2}, respectively. The corresponding Riesz fractional derivative and Riesz-Caputo derivative of (\ref{RI}) are defined by \cite{Agrawal1},
\[
{}^RD_t^\alpha f(t) = \frac{1}{\Gamma(n-\alpha)}\Big(\frac{d}{dt}\Big)^n \int_a^b |t - \tau|^{n-\alpha -1} f(\tau) \,d\tau, \; \; \; \alpha > 0,
\]
and
\[
{}^{C}D_t^\alpha f(t) = \frac{1}{\Gamma(n-\alpha)}\int_a^b |t - \tau|^{n-\alpha -1} \Big(\frac{d}{d\tau}\Big)^nf(\tau) \,d\tau, \; \; \; \alpha > 0,
\]
respectively. Further properties of these and similar operators can be found, for example, in \cite{key-2,key-8,key-9,hilferthreefold,Herrmann2,Torres2}.

This list is by no means complete and the interested reader may find detailed information about Davison and Essex \cite{D-S}, Coimbra \cite{Coimbra}, Marchaud \cite{Bapna,C-M,key-9}, Weyl-Marchaud \cite{W-M}, Chen-Marchaud \cite{C-M}, Marchaud-Hadamard \cite{hilferthreefold},  Miller-Ross \cite{key-11} and sequential fractional derivatives \cite{Furati,Loghmani,Chikrii} in the references given. 

Further details including comparison results of most of these derivatives can be found, for example, in the expository articles by Hilfer \cite{hilferthreefold}, and Atangana and Secer \cite{atan}. The historical notes about the developement of the theory can be found in the works of Ross \cite{Ross1, key-11,Ross2, Ross3, Ross4}. 

The other derivative that we elaborate in this paper, is the \emph{Hadamard fractional integral} introduced by J. Hadamard \cite{key-9,key-8,key-10}, and is given by,
\begin{equation}
\mathbf{I}^\alpha_{a+} f(x) = \frac{1}{\Gamma(\alpha)}\int_a^x \Bigg(\log\frac{x}{\tau}\Bigg)^{\alpha -1} f(\tau)\frac{d\tau}{\tau} ,
\label{eq:HI}
\end{equation}
for $Re(\alpha) > 0, \, x > a\geq 0$ while the \emph{Hadamard fractional derivative} of order $\alpha \in \mathbb{C}, Re(\alpha) >0$ is given by,
\begin{equation}
\mathbf{D}^\alpha_{a+} f(x) = \frac{1}{\Gamma(n-\alpha)}\Bigg(x\frac{d}{dx}\Bigg)^n \int_a^x \bigg(\log\frac{x}{\tau}\bigg)^{n-\alpha -1}f(\tau)\frac{d\tau}{\tau} ,
\label{eq:HD}
\end{equation}
for $x > a\geq 0$ where $n=\ipa$. The readers may have noticed that the version given in \cite[p. 111]{key-8} has misprints in its definitions of the Hadamard Derivatives. 

In 1940, generalizations of Riemann-Liouville fractional operators were introduced by Erd\'{e}lyi and Kober \cite{E-K}. \emph{Erd\'{e}lyi-Kober}-type fractional integral and derivative operators are defined by \cite{key-9,key-8,E-K,key-4,key-13,key-Y,Herrmann2},
\begin{align}
\;\;(I^\alpha_{a+; \,\rho, \,\eta} f)(x) &= \frac{\rho\,x^{-\rho(\alpha+\eta)}}{\Gamma(\alpha)}\int_a^x \frac{\tau^
{\rho\eta+\rho-1}\,f(\tau)}{(x^\rho-\tau^\rho)^{1-\alpha}}  d\tau ,  \label{eq:ek1}
\end{align} 
for $x > a\geq 0, Re(\alpha)> 0$ and 
\begin{align}
(D^\alpha_{a+; \,\rho, \,\eta} f)(x) = x^{-\rho\eta}\bigg(\frac{1}{\rho\,x^{\rho-1}}\,\frac{d}{dx}\bigg)^n x^{\rho(n+\eta)}\, \Big(I^{n-\alpha}_{a+; \,\rho, \,\eta+\alpha} f\Big)(x) ,
\label{eq:ek2}
\end{align}
for $x > a, Re(\alpha)\geq 0, \rho >0$. When $\rho =2, \, a=0$, the operators are called \emph{Erd\'{e}lyi-Kober} operators. When $\rho =1, \, a=0$, they are called \emph{Kober-Erd\'{e}lyi} or \emph{Kober operators} \cite[p.105]{key-8}. The geometric interpretation of this operator is discussed in \cite{Herrmann2}. 

In \cite{Osler2}, Osler further generalizes Erd\'{e}lyi's work by defining fractional integrals and derivatives of a function $f(x)$ with respect to another function $g(x)$ defined by \cite{key-9,key-8},
\begin{align}
\;\;(I^\alpha_{a+;g} f)(x) &= \frac{1}{\Gamma(\alpha)}\int_a^x \frac{g'(\tau)\,f(\tau)}{[(g(x) - g(\tau)]^{1-\alpha}}  d\tau , 
\label{eq:rd1}
\end{align}
for $x > a\geq 0, Re(\alpha)> 0$ and 
\begin{align}
(D^\alpha_{a+;g} f)(x) = \frac{1}{\Gamma(n-\alpha)}\bigg(\frac{1}{g'(x)}\,\frac{d}{dx}\bigg)^n \int_a^x \frac{g'(\tau)\,f(\tau)}{[(g(x) - g(\tau)]^{\alpha-n+1}}  d\tau ,
\label{eq:rd2}
\end{align}
for $x > a, Re(\alpha)\geq 0, \rho >0$, respectively. Where $g(x)$ is an increasing and positive function on $(a, \infty]$, having a continuos derivative $g(x)$ on $(a,\infty)$.

In \cite{key-0}, the author introduces a generalization to the Riemann-Liouville and Hadamard fractional integral and also provided existence results and semigroup properties. In the same reference, author introduces a generalised fractional derivative, which does not possess the inverse property. In this paper, we describe a new fractional derivative, which generalizes the Riemann-Liouville and the Hadamard fractional derivatives. Notice here that this new derivative possesses the inverse property [Theorem \,\ref{th:inv}].

\section{Generalization of the fractional integration and differentiation}\label{sec:2}
As in \cite{key-3}, consider the space $\textit{X}^p_c(a,b) \; (c\in \mathbb{R}, \, 1 \leq p \leq \infty)$ of those complex-valued Lebesgue measurable functions $f$ on $[a, b]$ for which $\|f\|_{\textit{X}^p_c} < \infty$, where the norm is defined by,
\begin{equation}\label{eq:dff1}
\|f\|_{\textit{X}^p_c} =\Bigg(\int^b_a |t^c f(t)|^p \frac{dt}{t}\Bigg)^{1/p} < \infty, 
\end{equation}
for $1 \leq p < \infty,\, c \in \mathbb{R}$ and for the case $p=\infty$,
\begin{equation} \label{eq:dff2}
\|f\|_{\textit{X}^\infty_c} = \text{ess sup}_{a \leq t \leq b} \big[t^c|f(t)|\big]  \quad ( c \in \mathbb{R}).
\end{equation}

We start with the definitions introduced in \cite{key-0} with a slight modification in the notation. Let $\Omega = [a,b], \, (-\infty <a < b < \infty)$ be a finite interval on the real axis, $\mathbb{R}$. The generalized fractional integral ${}^\rho I^\alpha_{a+}f$ of order $\alpha \in \mathbb{C} \; (Re(\alpha) > 0)$ of $f \in \textit{X}^p_c(a,b)$ is defined by,
\begin{equation}
\big({}^\rho \mathcal{I}^\alpha_{a+}f\big)(x) = \frac{\rho^{1- \alpha }}{\Gamma({\alpha})} \int^x_a \frac{\tau^{\rho-1} f(\tau) }{(x^\rho - \tau^\rho)^{1-\alpha}}\, d\tau, 
\label{eq:df1}
\end{equation}
for $x > a$, $Re(\alpha) > 0$, and $\rho > 0$. This integral is called the \emph{left-sided} fractional integral. Similarly, we can define the \emph{right-sided} fractional integral
${}^\rho I^\alpha_{b-}f$ by,
\begin{equation}
\big({}^\rho \mathcal{I}^\alpha_{b-}f\big)(x) = \frac{\rho^{1- \alpha }}{\Gamma({\alpha})} \int^b_x \frac{\tau^{\rho-1} f(\tau) }{(\tau^\rho - x^\rho)^{1-\alpha}}\, d\tau,
\label{eq:df2}
\end{equation} 
for $x < b$ and $Re(\alpha) > 0$. These are the fractional generalizations of the $n-$fold left- and right- integrals of the form $ \int_a^x t_1^{\rho-1} dt_1 \int_a^{t_1} t_2^{\rho-1} dt_2 \cdots \int_a^{t_{n -1}} t_n^{\rho-1} f(t_n)dt_n$ and  $\int_x^b t_1^{\rho-1} dt_1 \int^b_{t_1} t_2^{\rho-1} dt_2 \cdots \int^b_{t_{n -1}} t_n^{\rho-1} f(t_n)dt_n $ for $n \in \mathbb{N}$, respectively. When $b=\infty$, the generalized fractional integral is called a Liouville-type integral and the case $a=\infty$ is referred to as the Weyl-derivative \cite{weyl1}.

\begin{remark}
It can be seen that these operators are not equivalent to Erd\'{e}lyi-Kober operators or equation $(\ref{eq:rd1})$ when $g(x)=x^\rho$, but is different from a factor of $\rho^{-\alpha}$, which is essential in the case of the Hadamard operators. 
\end{remark}

It is worthy to mention here that we interchangeably use the notations ${}^\rho_a \mathcal{I}^\alpha_x$ and ${}^\rho \mathcal{I}^\alpha_{a+}$ for the geneneralized integral (\ref{eq:df1}). We obtained conditions for the integration operator ${}^\rho_a \mathcal{I}^\alpha_x$ to be bounded in $\textit{X}^p_c(a,b)$, and also established semigroup property for the generalized fractional integration operator. Those two results are summarized bellow \cite{key-0}.
\begin{theorem}\label{eq:th1}
Let $\alpha > 0,\, 1 \leq p \leq \infty,\, 0 <a < b < \infty$ and let $\rho \in \mathbb{R}$ and $c \in \mathbb{R}$ be such that $\rho -1 \geq c$. Then the operator ${}^\rho_a I^\alpha_t$ is bounded in $\textit{X}^p_c(a,b)$ and 
\begin{equation}
\|{}^\rho_a \mathcal{I}^\alpha_{t}f \|_{\textit{X}^p_c} \leq \mathcal{K}\|f\|_{\textit{X}^p_c}, \nonumber
\end{equation}
\label{eq:thm1}
where 
\begin{equation}
\mathcal{K} = \frac{b^{\alpha\rho-1}}{\Gamma(\alpha)}\int^{\frac{b}{a}}_1 u^{c-\alpha\rho-1}\Bigg(\frac{u^{\rho}-1}{\rho}\Bigg)^{\alpha -1} du, \quad \text{$\rho \neq 0.$} \nonumber
\end{equation}
\label{eq:con1}
\end{theorem}
Kilbas proved a version of Theorem \ref{eq:th1} in \cite{key-3}, for the special case when $\rho \rightarrow 0^+$. In \cite{key-0}, we proved the boundedness of the operator ${}^\rho_a \mathcal{I}^\alpha_{t}$ in the space $\textit{L}^p (a,b)$, and the Semigroup property, \textsl{i.e.}, for $\alpha >0,\, \beta >0,\, 1 \leq p \leq \infty, \, 0 < a < \infty, \rho \in \mathbb{R}, c \in \mathbb{R}$ and $\rho \geq c$, ${}^\rho_{a} \mathcal{I}^{\alpha}_{t}\,{}^\rho_{a} \mathcal{I}^{\beta}_{t}f = {}^\rho_{a} \mathcal{I}^{\alpha+\beta}_{t}f$ for $f \in \textit{X}^p_c(a,b)$, was also obtained. The interested reader is refered to \cite{key-0} for further results related to the generalized fractional integral in question.

Next we give the main result of this paper. First consider the generalized fractional derivatives defined below.
\begin{definition}(Generalized Fractional Derivatives)\\
Let $\alpha \in \mathbb{C}, Re(\alpha) \geq 0, n=\ipa$ and $\rho >0.$ The generalized fractional derivatives, corresponding to the generalized fractional integrals (\ref{eq:df1}) and (\ref{eq:df2}), are defined, for $0 \leq a < x < b \leq \infty$, by
\begin{align}
\big({}^\rho \mathcal{D}^\alpha_{a+}f\big)(x)&= \bigg(x^{1-\rho} \,\frac{d}{dx}\bigg)^n\,\, \big({}^\rho \mathcal{I}^{n-\alpha}_{a+}f\big)(x)\nonumber\\
 &= \frac{\rho^{\alpha-n+1 }}{\Gamma({n-\alpha})} \, \bigg(x^{1-\rho} \,\frac{d}{dx}\bigg)^n \int^x_a \frac{\tau^{\rho-1} f(\tau) }{(x^\rho - \tau^\rho)^{\alpha-n+1}}\, d\tau,
\label{eq:gd1}
\end{align}
and
\begin{align}
\big({}^\rho \mathcal{D}^\alpha_{b-}f\big)(x) &= \bigg(-x^{1-\rho} \,\frac{d}{dx}\bigg)^n\,\, \big({}^\rho \mathcal{I}^{n-\alpha}_{b-}f\big)(x)\nonumber\\
 &= \frac{\rho^{\alpha-n+1 }}{\Gamma({n-\alpha})}\bigg(-x^{1-\rho}\frac{d}{dx}\bigg)^n \int^b_x\frac{\tau^{\rho-1} f(\tau) }{(\tau^\rho - x^\rho)^{\alpha-n+1}}\, d\tau,
\label{eq:gd2}
\end{align}
if the integrals exist. 
\end{definition}




In the next section we give several properties pertaining to the generalized fractional derivative we have just defined. In Theorem \ref{th1}, we show that the Riemann-Liouville and the Hadamard fractional derivatives are special cases of the generalized derivative in question.

\section{Properties of the generalized operators}\label{sec:3}
In this section we introduce several identites pertaining to the generalized fractional operators. Mainly, we give the inverse property, composition theorem, index theorem and linearity followed by Taylor-like expansion, which extends the classical Taylor series to include fractional derivatives. 

The theorem below gives the relations of generalized fractional derivatives to that of Riemann-Liouville and Hadamard. For simplicity we give only the \emph{left-sided} versions here. 
\begin{theorem}\label{th1}
Let $\alpha \in \mathbb{C},\, Re(\alpha) \geq 0,\, n=\lceil Re(\alpha) \rceil$ and $\rho >0.$ Then, for $x > a$,
\begin{align}
\noindent	&1. \lim_{\rho \rightarrow 1}\big({}^\rho \mathcal{I}^\alpha_{a+}f\big)(x) =\frac{1}{\Gamma(\alpha)}\int_a^x (x - \tau)^{\alpha -1} f(\tau) d\tau, \label{eq:th5} & & & &\\
	&2. \lim_{\rho \rightarrow 0^+}\big({}^\rho \mathcal{I}^\alpha_{a+}f\big)(x)=\frac{1}{\Gamma(\alpha)}\int_a^x \Bigg(\log\frac{x}{\tau}\Bigg)^{\alpha -1} f(\tau)\frac{d\tau}{\tau}, \label{eq:th2} & & & &\\
	&3. \lim_{\rho \rightarrow 1}\big({}^\rho \mathcal{D}^\alpha_{a+}f\big)(x)=\bigg(\frac{d}{dx}\bigg)^n \,\frac{1}{\Gamma({n-\alpha})} \, \int^x_a \frac{ f(\tau) }{(x - \tau)^{\alpha-n+1}}\, d\tau,\label{eq:th3} & & & &\\
	&4. \lim_{\rho \rightarrow 0^+}\big({}^\rho \mathcal{D}^\alpha_{a+}f\big)(x)=\frac{1}{\Gamma(n-\alpha)}\Bigg(x\frac{d}{dx}\Bigg)^n \int_a^x \Bigg(\log\frac{x}{\tau}\Bigg)^{n-\alpha -1} f(\tau)\frac{d\tau}{\tau}.\label{eq:th4}& & & &
\end{align}
\end{theorem}
\begin{proof}
The equations (\ref{eq:th5}) and (\ref{eq:th3}) follow from taking limits when $\rho \rightarrow 1$, while (\ref{eq:th2}) follows from the L'Hospital rule by noticing that \cite{key-0}
\begin{align*}
\lim_{\rho \rightarrow 0^+}\frac{\rho^{1- \alpha }}{\Gamma({\alpha})} \int^x_a &\frac{ f(\tau)\tau^{\rho-1} }{(x^\rho - \tau^\rho)^{1-\alpha}}\, d\tau \\
  &= \frac{1}{\Gamma({\alpha})}\int^x_a \lim_{\rho \rightarrow 0^+} f(\tau)\tau^{\rho-1} \Bigg(\frac{x^\rho - \tau^\rho}{\rho}\Bigg)^{\alpha -1}d\tau , \\
	&= \frac{1}{\Gamma(\alpha)}\int_a^x \Bigg(\log\frac{x}{\tau}\Bigg)^{\alpha -1} f(\tau)\frac{d\tau}{\tau} .
\end{align*}
The proof of (\ref{eq:th4}) is similar.
\end{proof}
Similar results for \emph{right-sided} integrals and derivatives also exist and can be proved similarly. 
Note that the equations (\ref{eq:th5}) and (\ref{eq:th3}) are related to the Riemann-Liouvile operators, while equaitons (\ref{eq:th2}) and (\ref{eq:th4}) are related to the Hadamard operators. The results similar to equations (\ref{eq:th5}) and (\ref{eq:th3}) can also be derived from the Erd\'{e}lyi-Kober operators, though it is not possible to obtain equivalent results for (\ref{eq:th2}) and (\ref{eq:th4}). This is due to the absence of the factor $\rho^{\alpha}$, which is needed in the limit to obtain the Hadamard-type operators.

We are now ready to state and prove the following properties. The first is the inverse property. 
\begin{theorem}[Inverse property]\label{th:inv}
Let $0 <\alpha <1$, and $f \in \textit{X}^p_c(a,b), \,\, \rho >0$. Then, for $a>0,\, \rho >0$,
\begin{equation}
  \Big({}^\rho \mathcal{D}^\alpha_{a+}\,{}^\rho \mathcal{I}^\alpha_{a+}\Big)f(x) = f(x).
\end{equation}
\end{theorem}
\begin{proof}

We prove this using Fubini's theorem and Dirichlet technique \cite[p. 59]{key-2}. From direct integration, we have
\begin{align*}
\Big({}^\rho \mathcal{D}^\alpha_{a+}\,&{}^\rho \mathcal{I}^\alpha_{a+}\Big)f(x) \\
   &=\frac{\rho^\alpha}{\Gamma({1-\alpha})} \, \bigg(x^{1-\rho} \,\frac{d}{dx}\bigg) \int^x_a \frac{\tau^{\rho-1}}{(x^\rho - \tau^\rho)^\alpha}\,\cdot \frac{\rho^{1- \alpha }}{\Gamma({\alpha})} \int^\tau_a \frac{s^{\rho-1} f(\tau) }{(\tau^\rho - s^\rho)^{1-\alpha}}\, ds\, d\tau, \\
   &= \frac{\rho}{\Gamma({1-\alpha})\Gamma({\alpha})} \, \bigg(x^{1-\rho} \,\frac{d}{dx}\bigg) \int^x_a f(s)\,s^{\rho-1} \int^x_s \frac{(\tau^\rho - s^\rho)^{\alpha-1}}{(x^\rho - \tau^\rho)^\alpha} \tau^{\rho-1} \, d\tau ds, \\
   &= \frac{\rho}{\Gamma({1-\alpha})\Gamma({\alpha})} \, \bigg(x^{1-\rho} \,\frac{d}{dx}\bigg) \int^x_a f(s)\,s^{\rho-1}ds \cdot \frac{\Gamma({1-\alpha})\Gamma({\alpha})}{\rho},\\
   &= f(x).\nonumber
\end{align*}
Notice here that the inner integral is evaluated by the change of variable, $t=(\tau^\rho-s^\rho)/(x^\rho-s^\rho)$, and using the Beta function defined by, $B(\alpha,\beta)=\int^1_0 t^{\alpha-1}(1-t)^{\beta-1} dt$ and the fact that $B(\alpha,\beta)=\Gamma(\alpha)\Gamma(\beta)/\Gamma(\alpha+\beta)$. This completes the proof.
\end{proof}


Next is the index theorem. 
\begin{theorem}[Index theorem]
Let $\alpha,\, \beta \in \mathbb{C}$ be such that $0 < Re(\alpha)<1$ and $0< Re(\beta) < 1.$ If $0 < a < b < \infty$ and $ 1 \leq p \leq \infty$, then, for $f \in \textit{X}^p_c(a,b), \,\, \rho >0$,
\begin{equation}
   {}^\rho \mathcal{I}^\alpha_{a+}\,{}^\rho \mathcal{I}^\beta_{a+}\,f = {}^\rho \mathcal{I}^{\alpha + \beta}_{a+}\,f \;\;\text{and}\;\;{}^\rho \mathcal{D}^\alpha_{a+}\,{}^\rho \mathcal{D}^\beta_{a+}\,f = {}^\rho \mathcal{D}^{\alpha + \beta}_{a+}\,f. \nonumber
\end{equation}
\end{theorem}
\begin{proof}
The interested reader may find the proof of the first result in \cite{key-0}. The proof of the second result follows from direct integration.
\end{proof}

Now we turn our attention to the following. The compositions between the generalized fractional differentiation and the generalized fractional integration are given by the following result.
\begin{theorem}[Composition]
Let $\alpha,\, \beta \in \mathbb{C}$ be such that $0 < Re(\alpha) < Re(\beta) < 1.$ If $0 < a < b < \infty$ and $ 1 \leq p \leq \infty$, then, for $f \in L^p(a,b), \,\, \rho >0$,
\begin{equation}
  {}^\rho \mathcal{D}^\alpha_{a+}\,{}^\rho \mathcal{I}^\beta_{a+}\,f = {}^\rho \mathcal{I}^{\beta - \alpha}_{a+}\,f \;\;\text{and}\;\;{}^\rho \mathcal{D}^\alpha_{b-}\,{}^\rho \mathcal{I}^\beta_{b-}\,f = {}^\rho \mathcal{I}^{\beta - \alpha}_{b-}\,f. \nonumber
\end{equation}
\end{theorem}

\begin{proof}
The proof is similar to that of Theorem \ref{th:inv}. Notice that we use the property $\Gamma(x+1) = x\,\Gamma(x)$ of the Gamma function here. 
\end{proof}
\begin{remark} 
The author suggests the interested readers refer Property 2.27 in \cite{key-8} for similar properties of the Hadamard fractional operator. 
\end{remark}

As in the case of classical derivatives, the generalized operators also satisfy the linearity property. 
\begin{theorem}[Linearity property]
Let $\alpha \in \mathbb{C}$ be such that $0 < Re(\alpha)<1$. If $0 < a < b < \infty$ and $ 1 \leq p \leq \infty$, then for $f, g \in \textit{X}^p_c(a,b)$ and $\rho >0$,
\begin{equation}
   {}^\rho \mathcal{I}^\alpha_{a+}\big(f+g\big) = {}^\rho \mathcal{I}^\alpha_{a+}\,f + {}^\rho \mathcal{I}^\alpha_{a+}\,g
   \label{lin-1}
\end{equation}
and,
\begin{equation}
{}^\rho \mathcal{D}^\alpha_{a+} \big(f+g\big) = {}^\rho \mathcal{D}^\alpha_{a+}\,f + {}^\rho \mathcal{D}^\alpha_{a+}\,g. 
\label{lin-2}
\end{equation}
\end{theorem}
\begin{proof}
The results follow from direct integration. 
\end{proof}

\section{Taylor-like expansions}\label{sec:4}
The identities (\ref{lin-1}) and (\ref{lin-2}) can be used to express the generalized fractional derivative of an analytic function, which possesses a Taylor series expansion.
\begin{theorem}
Let $\mathcal{R}$ be a simply connected region containing the origin. Also, let $f(z)$ be analytic in $\mathcal{R}$ and $\rho >0$. 
Then,
\begin{align}
{}^\rho D^\alpha_{0+} f(z) &= \frac{\rho^{\alpha-1}}{z^{\alpha\rho}}\sum^\infty_{i=0}\frac{\Gamma\Big(1+\frac{i}{\rho}\Big)\,}{i!\, \Gamma\Big(1+\frac{i}{\rho}-\alpha\Big)}f^{(i)}(0)\, z^{i}, \quad \text{for} \; z\neq 0.
\label{taylor-thm} 
\end{align}
\end{theorem}
\begin{proof}
The proof immediately follows from (\ref{eq:gen}) using the Taylor expansion of $f(z)$ near $0$.
\end{proof}
\begin{remark}
One of the advantages of this result is that it enables us to find the generalized fractional derivative of any function, which is analytic in a neighborhood of $0$. The direct calculation of the derivative may be too much involved, for example consider the case of $f(x)=\sin^{-1}x.$  
\end{remark}

In \cite{u-1}, Gaboury et al. obtained two representations of ${}^\rho \mathcal{I}^\alpha_{0+}$ given by the following result.
\begin{theorem} Let $\mathcal{R}$ be a simply connected region containing the origin. Also, let $f(z)$ be analytic in $\mathcal{R}$. Then, for $\alpha\in \mathbb{C}$ and $\rho \in \mathbb{C}$ with $Re(\alpha)>0$ and $Re(\rho)>0$, the following relations hold true
\begin{align*}
    {}^\rho \mathcal{I}^\alpha_{0+}f(z) = \frac{\rho^{1-\alpha}}{\Gamma(\alpha)}z^{\rho(\alpha-1)}&\sum^\infty_{n=0}\frac{(1-\alpha)_n}{n!}z^{-n\rho}\\
		         &\times \Gamma\big(\rho(n+1)\big)D_z^{-\rho(n+1)}f(w-z)\Big|_{w=z},
 \end{align*}	
and
\begin{align*}
    {}^\rho \mathcal{I}^\alpha_{0+}f(z) = \frac{\rho^{1-\alpha}}{\Gamma(\alpha)}z^{\rho(\alpha-1)}&\sum^\infty_{n=0}\frac{(1-\alpha)_n}{n!}z^{-n\rho}
		         D_z^{-1}z^{\rho(n+1)-1}f(z),
 \end{align*}	
where $(\lambda)_n$ is the Pochhammer's symbol given by $(\lambda)_n = \Gamma(\lambda+n)/\Gamma(\lambda), \; (\lambda)_0=1$.
\end{theorem}			

In this paper we obtain similar results for the fractional derivative operator ${}^\rho \mathcal{D}^\alpha_{0+}$ given by the next result. 

\begin{theorem} Let $\mathcal{R}$ be a simply connected region containing the origin. Also, let $f(z)$ be analytic in $\mathcal{R}$. Then, for $\alpha \in \mathbb{C}$ with $0< Re(\alpha)<1$ and $\rho \in \mathbb{C}$ with $Re(\rho)>0$, the following relations hold true.
\begin{align}
    {}^\rho \mathcal{D}^\alpha_{0+}f(z) = \frac{\rho^{\alpha}}{\Gamma(1-\alpha)}\Big(z^{1-\rho}\frac{d}{dz}\Big)&\sum^\infty_{k=0}\frac{(\alpha)_k}{k!}z^{-(\alpha + k)\rho} \nonumber\\ 
		         &\times \Gamma\big(\rho(k+1)\big)I_{0+}^{\rho(k+1)}f(w-z)\Big|_{w=z}, \label{th2-1}
 \end{align}	
and
\begin{align}
    {}^\rho \mathcal{D}^\alpha_{0+}f(z) &= \frac{\rho^{\alpha}}{\Gamma(1-\alpha)}\sum^\infty_{k=0}\frac{(\alpha)_k}{k!}\frac{\Gamma(\rho(k+1))}{z^{\rho(\alpha +k+1)}} \nonumber\\
		         &\times \Bigg\{zD_{0+}^{1-\rho(k+1)}f(w-z)\Big|_{w=z} - \rho(k+\alpha)I_{0+}^{\rho(k+1)}f(w-z)\Big|_{w=z}\Bigg\},\label{th2-2}
 \end{align}	
where $(\lambda)_n$ is the Pochhammer's symbol given by $(\lambda)_n = \Gamma(\lambda+n)/\Gamma(\lambda), \; (\lambda)_0=1$ and $I_{0+}^{\alpha}$ and $D_{0+}^{\alpha}$ are the Riemann-Liouville fractional integral and derivative of variable $z$, respectively. 
\end{theorem}		

\begin{proof} For the proof, we follow a similar approach as appeared in \cite[page 4]{u-1}. Consider the generalized operator ${}^\rho \mathcal{D}^\alpha_{0+}$ in the complex plane with $0< Re(\alpha)<1$. Then $n=\lceil Re(\alpha) \rceil = 1$ and by (\ref{eq:gd1}),
\begin{align*}
\big({}^\rho \mathcal{D}^\alpha_{a+}f\big)(z)&= \frac{\rho^{\alpha}}{\Gamma({1-\alpha})} \, \bigg(z^{1-\rho} \,\frac{d}{dz}\bigg) \int^z_0 \frac{\tau^{\rho-1} f(\tau) }{(z^\rho - \tau^\rho)^{\alpha}}\, d\tau
\end{align*}	
Making the change of variable $\tau=z-\xi, \, (d\tau=-d\xi)$, we have	
\begin{align}
\big({}^\rho \mathcal{D}^\alpha_{a+}f\big)(z)&= \frac{\rho^{\alpha}}{\Gamma({1-\alpha})} \, \bigg(z^{1-\rho} \,\frac{d}{dz}\bigg) \int^z_0 \frac{f(z-\xi)(z-\xi)^{\rho-1}}{\big(z^\rho - (z-\xi)^\rho\big)^{\alpha}}\, d\xi, \nonumber\\
  &= \frac{\rho^{\alpha}}{\Gamma({1-\alpha})} \, \bigg(z^{1-\rho} \,\frac{d}{dz}\bigg) z^{-\rho\alpha} \int^z_0 \frac{f(z-\xi)(z-\xi)^{\rho-1}  }{\Big(1-\big(\frac{z-\xi}{z}\big)^\rho\Big)^\alpha}\, d\xi, \nonumber\\
  &= \frac{\rho^{\alpha}}{\Gamma({1-\alpha})} \, \bigg(z^{1-\rho} \,\frac{d}{dz}\bigg) z^{-\rho \alpha} \int^z_0 \sum^\infty_{k=0}\frac{(\alpha)_k}{k!}\bigg(\frac{z-\xi}{z}\bigg)^{\rho k} \frac{f(z-\xi)}{(z-\xi)^{1-\rho}} \, d\xi, \label{eq-1}\\
  &= \frac{\rho^{\alpha}}{\Gamma({1-\alpha})} \, \bigg(z^{1-\rho} \,\frac{d}{dz}\bigg) \sum^\infty_{k=0}\frac{(\alpha)_k}{k!} z^{-\rho(\alpha +k)} \int^z_0 f(z-\xi)(z-\xi)^{\rho(k+1)-1}\, d\xi. \label{eq-2}
\end{align}
Rewriting the integral in (\ref{eq-2}) in terms of the Riemann-Liouville fractional integral operator yields (\ref{th2-1}). Notice that in (\ref{eq-1}), we have used the power series expansion of $\big(1-\big(\frac{z-\xi}{z}\big)^\rho\big)^\alpha$ 	for the values of $\xi$ near $z$. Equation (\ref{th2-2}) follows from (\ref{eq-1}) after applying the product rule and rewriting the resulted sum in terms of the Riemann-Liouville integral and derivative. 
\end{proof}	
		
\begin{figure}[h]
  \title{	The Generalized Derivatives of the Power Function $f(x)=x^\nu$}	
	\centering
	  \subfloat[$\nu$= 1.0]{\includegraphics[width=2.3in, height=2.2in,  trim=0 0 0 0]{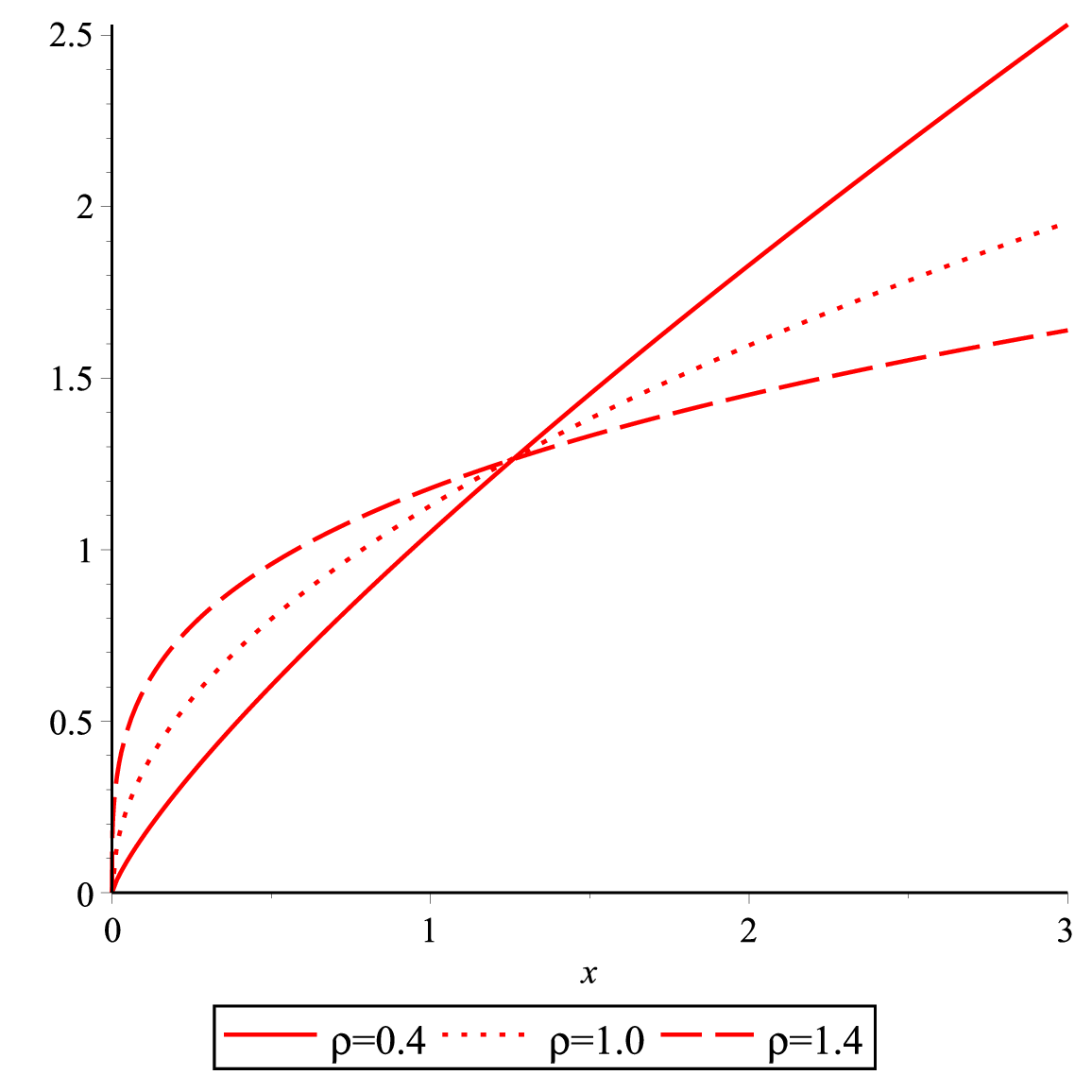}} 
	  \subfloat[$\nu$= 2.0]{\includegraphics[width=2.3in, height=2.2in,  trim=0 0 0 0]{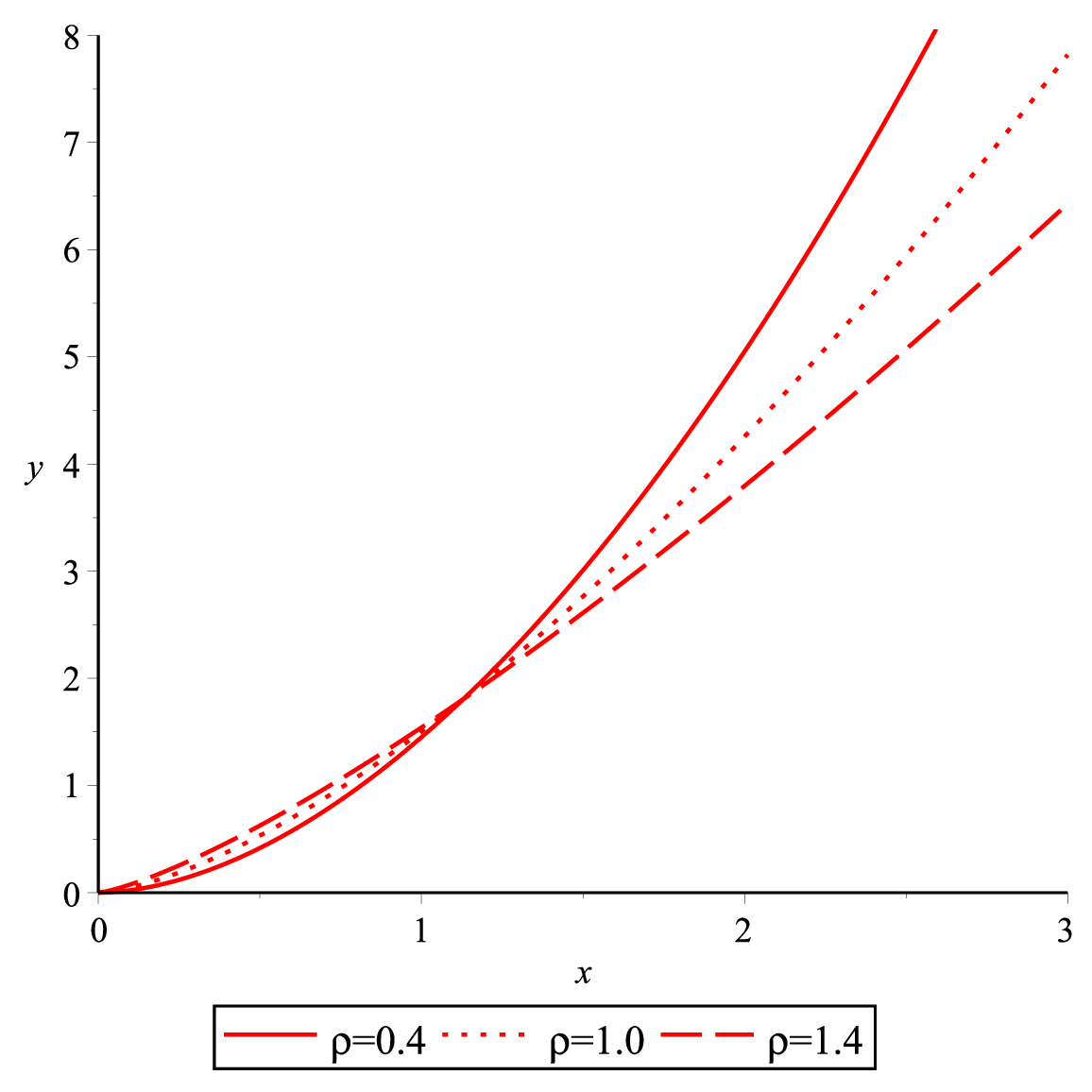}}\\
	  
	  \subfloat[$\nu$= 0.5]{\includegraphics[width=2.3in, height=2.2in,  trim = 0 0 0 0]{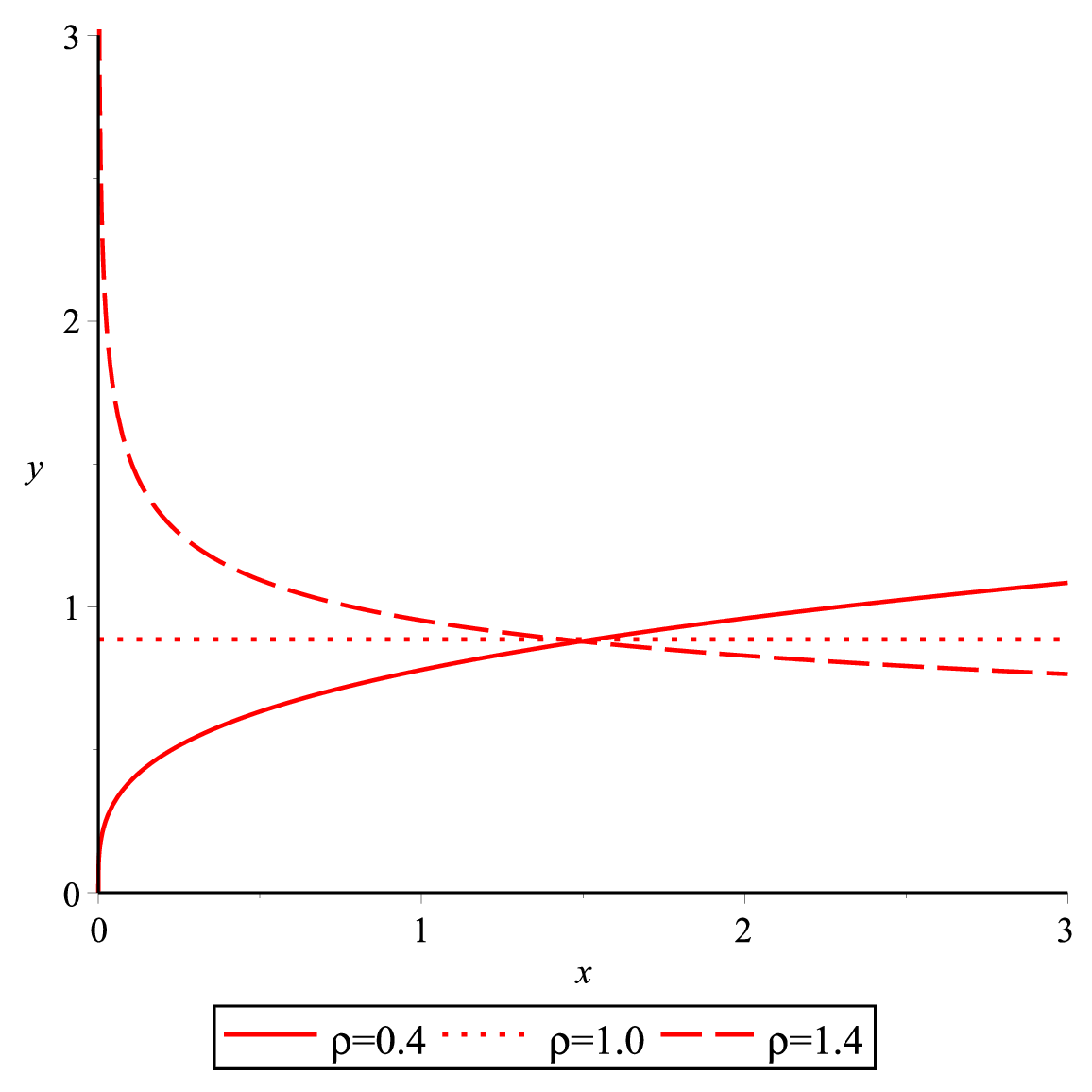}}
	  \subfloat[$\nu$= 1.5]{\includegraphics[width=2.3in, height=2.2in,  trim = 0 0 0 0]{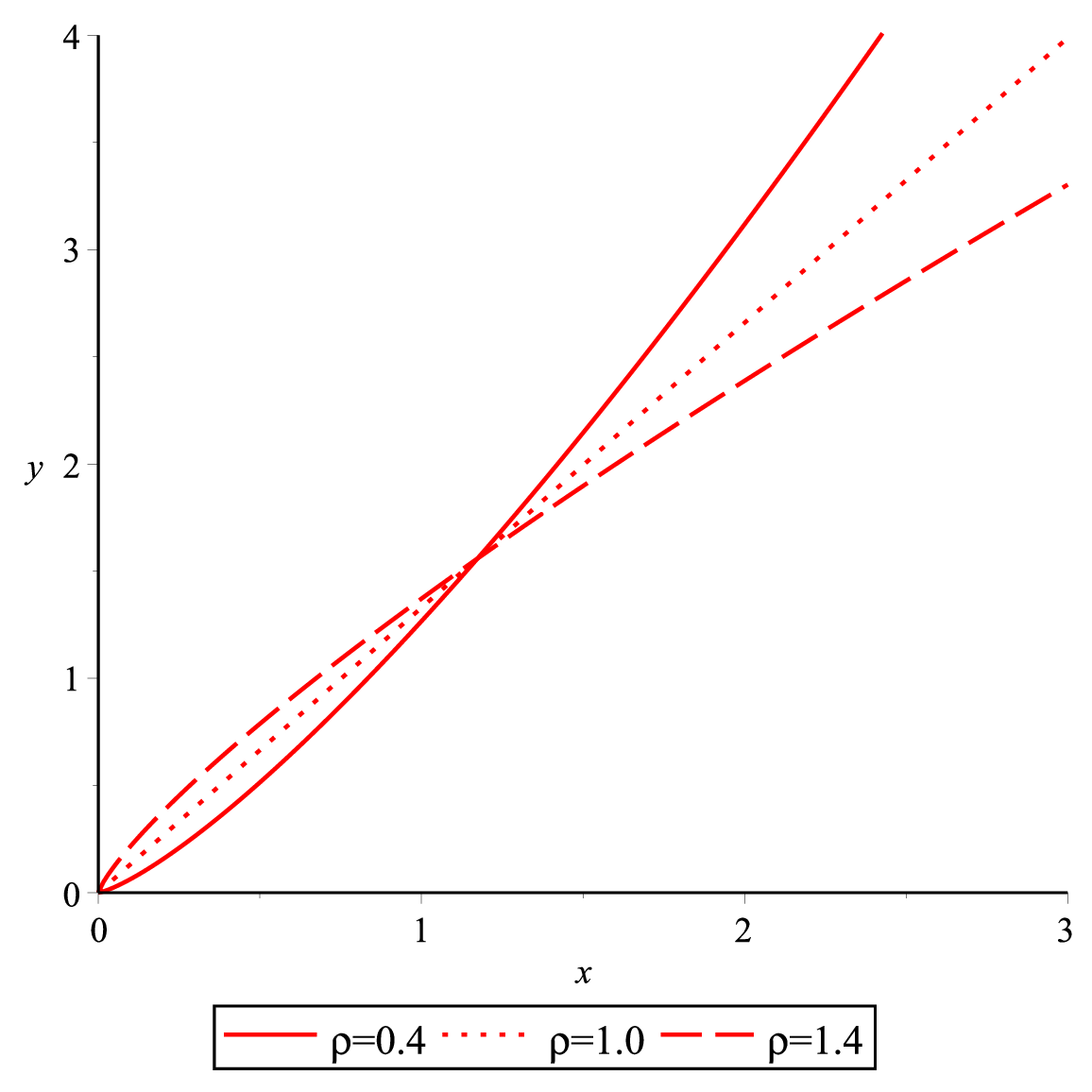}}
	\caption{Generalized fractional derivative of the power function $f(x)=x^\nu$ for $\rho = 0.4,\, 1.0,\, 1.4$ and $\nu = 0.5,\, 1.0, \, 1.5, \, 2.0.$}  
	\label{fig:FD-1}
\end{figure}
Let us consider an example to illustrate the results. We shall find the generalized fractional derivative (\ref{eq:gd1}) of the power function and investigate the behavior for different values of $\rho, \alpha$ and $\nu$. For simplicity assume $\alpha \in \mathbb{R}^+, \, 0 < \alpha < 1$ and $a=0$. 
\begin{example}
We consider  the function $f(x) = x^\nu$, where $\nu \in \mathbb{R}$. The formula (\ref{eq:gd1}) yields
\begin{align}
  {}^\rho D^\alpha_{0+} x^\nu = \frac{\rho^\alpha}{\Gamma(1-\alpha)} \, \Bigg(x^{1-\rho} \frac{d}{dx}\Bigg)\int^x_0 \frac{t^{\rho-1}}{(x^{\rho}-t^{\rho})^\alpha}\, t^\nu \, dt.
\label{eq:simple1}  
\end{align}
To evaluate the inner integral, use the substitution $u=t^{\rho}/x^{\rho}$ to obtain,
\begin{align}
   \int^x_0 \frac{t^{\rho-1}}{(x^{\rho}-t^{\rho})^\alpha}\, t^\nu \, dt 
       &= \frac{x^{\nu +\rho(1-\alpha)}}{\rho} \int^1_0 \frac{u^\frac{\nu}{\rho}}{(1-u)^\alpha} \, du, \nonumber \\
       &= \frac{x^{\nu +\rho(1-\alpha)}}{\rho}\, B\Big(1-\alpha, 1+\frac{\nu}{\rho}\Big),\nonumber
\end{align}
where $B(.\, ,.)$ is the Beta function. Thus, we obtain,
\begin{align}
    {}^\rho D^\alpha_{0+} x^\nu &= \frac{\Gamma\Big(1+\frac{\nu}{\rho}\Big)\,\rho^{\alpha-1}}{\Gamma\Big(1+\frac{\nu}{\rho}-\alpha\Big)}\, x^{\nu -\alpha\rho}, 
\label{eq:gen}    
\end{align}
for $\rho > 0$, after using the properties of the Beta function \cite{key-8} and the relation $\Gamma(z+1)=z\,\Gamma(z)$. When $\rho = 1$ we obtain the Riemann-Liouville fractional derivative of the power function given by \cite{key-9,key-8,key-2},
\begin{equation}
  {}^1D^\alpha_{0+} x^\nu = \frac{\Gamma\Big(1+\nu\Big)}{\Gamma\Big(1+\nu-\alpha\Big)}\; x^{\nu-\alpha}. \label{eq:eg1}
\end{equation}
\end{example}
Equation (\ref{eq:eg1}) agrees well with the standard results obtained for Riemann-Liouville fractional derivative (\ref{eq:lRLD}). Interestingly enough, for $\alpha = 1,\, \rho = 1$, we obtain ${}^1D^1_{0+} x^\nu = \nu\,x^{\nu -1}$, as one would expect.

\begin{figure}[h]
	\centering
	  \subfloat[$\nu$= 2.0]{\includegraphics[width=2.3in, height=2.2in,  clip=true, trim = 0 0 0 0]{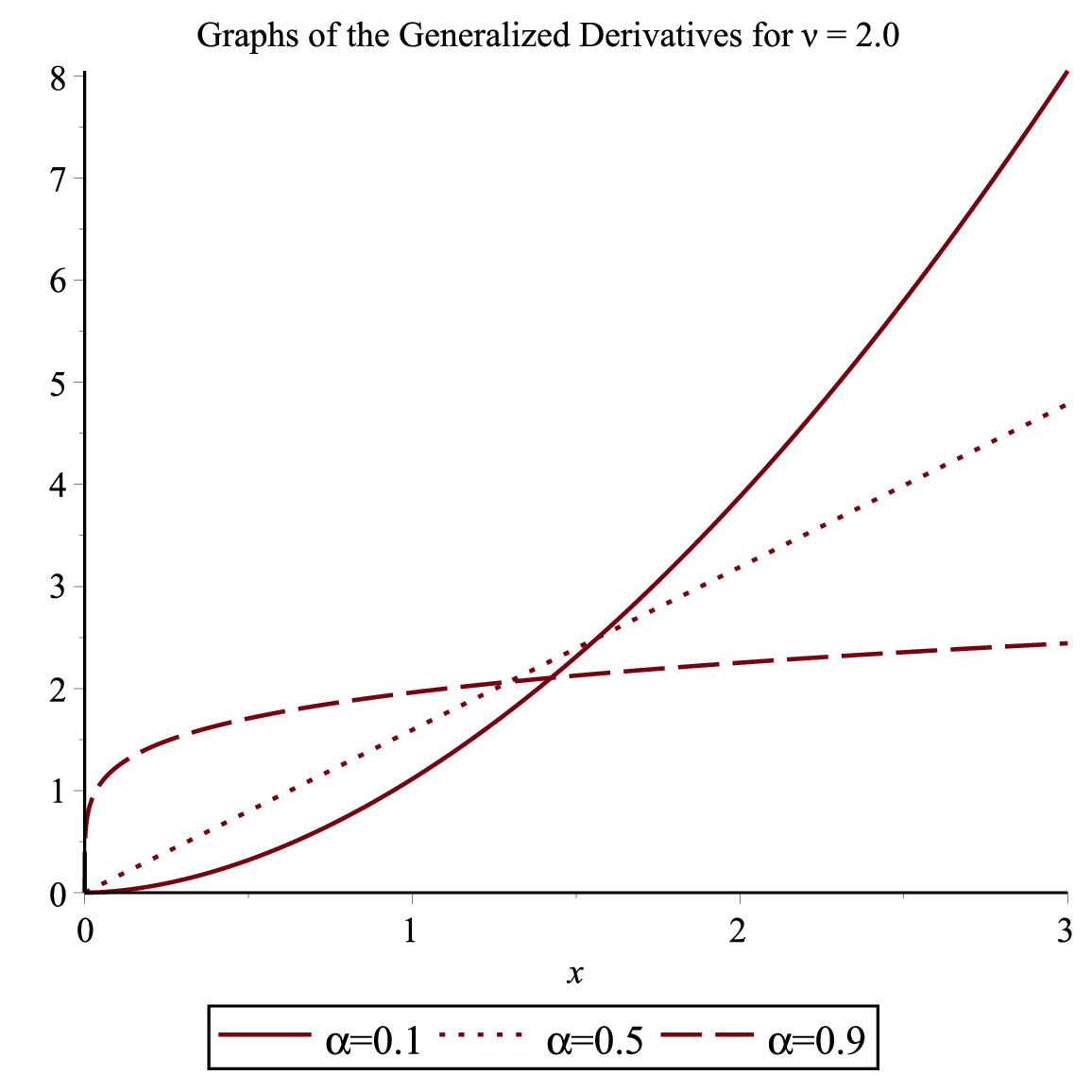}}
	  \subfloat[$\alpha$= 0.5]{\includegraphics[width=2.3in, height=2.2in,clip=true, trim = 0 0 0 0]{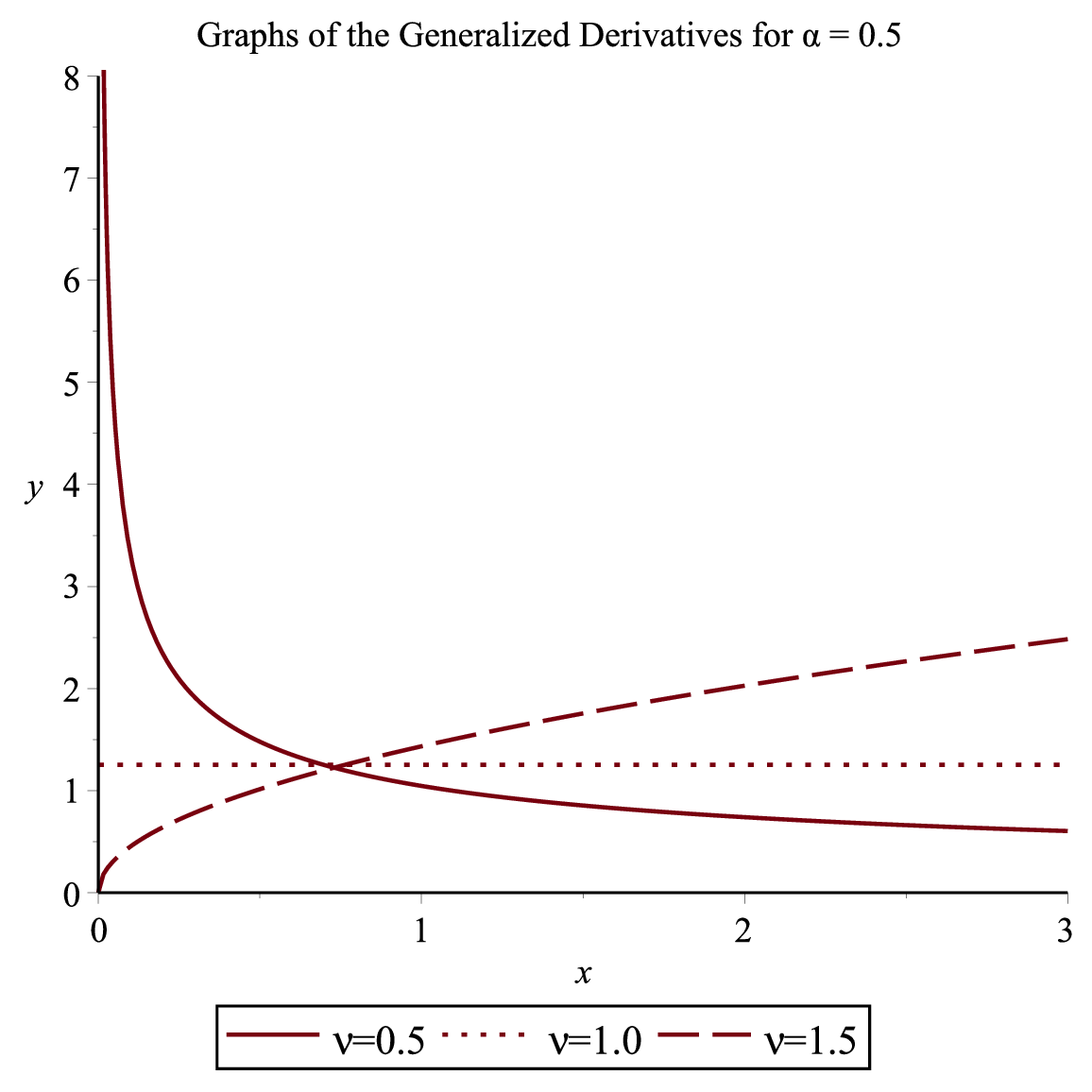}}  
	\caption{Generalized fractional derivative of the power function $f(x)=x^\nu$ for $\alpha = 0.1,\, 0.5,\, 0.9$ and $\nu = 0.5, \, 1.0, \, 1.5$}        
	\label{fig:FD-2}
\end{figure}

To compare the results, we plot (\ref{eq:gen}) for several values of $\rho\in \mathbb{R}$. We also consider different values of $\nu$ to see the effect on the degree of the power function. 
Figure \ref{fig:FD-1} summaries the comparision results for $\rho$ and $\nu$, while Figure \ref{fig:FD-2} summaries the comparision results for different values of $\alpha$ and $\nu$. We notice that the characteristics of the fractional derivative are mainly depend on the values of $\rho$, thus it provides a new direction for applications. In Figure 1(c) and 1(d), we further notice that for $\nu=0.5$ and $\nu=1.5$, the concavity of the function changes when $\rho$ changes near 1, at which the fractional derivative takes a constant value for different values of $x$ or becomes linear. 

\begin{remark}
It is easy to see that the generalized fractional integrals and derivatives in (\ref{eq:df1}), (\ref{eq:df2}), (\ref{eq:gd1}) and (\ref{eq:gd2}) can be expressed in terms of Erd\'{e}lyi-Kober operator given in (\ref{eq:ek1}) and (\ref{eq:ek2}). According to Theorem 4.1, it is evident that the Riemann-Liouville and the Hadamard fractional derivatives are special cases of this new derivative. It is not possible to obtain the Hadamard operators in the case of the Erd\'{e}lyi-Kober operators due to a reason explained earlier. 
\end{remark}




\section{Conclusion}

The paper presents an extended fractional differentiation, which generalizes the Riemann-Liouville and the Hadamard fractional derivatives into a single form, which when a parameter is fixed at different values or by taking limits produces the above derivatives as special cases. It is also pointed out that this new derivative is not equivalent to Erd\'{e}lyi-Kober derivative, which is the genealization of the Riemann-Liouville derivative. 

\section*{Acknowledgements}
  The author would like to thank Jerzy Kocik, the Department of Mathematics at Southern Illinois University-Carbondale, for his valuable comments and suggestions. 
	


\bibliographystyle{abbrv}
\bibliography{refs}      













\end{document}